\documentclass{article}
\usepackage[utf8]{inputenc}
\usepackage{amsmath}
\usepackage{amsfonts}
\usepackage{dirtytalk}

\title{A simple mnemonic to compute sums of powers}
\author{A. Mariani}

\begin{document}

\maketitle

\begin{abstract}
We give a simple recursive formula to obtain the general sum of the first $N$ natural numbers to the $r$th power. Our method allows one to obtain the general formula for the $(r+1)$th power once one knows the general formula for the $r$th power. The method is very simple to remember owing to an analogy with differentiation and integration. Unlike previously known methods, no knowledge of additional specific constants (such as the Bernoulli numbers) is needed. This makes it particularly suitable for applications in cases when one cannot consult external references, for example mathematics competitions.

\end{abstract}

\section{Introduction}

Sums of powers have fascinated mathematicians for centuries. The sum of the first $N$ natural numbers is given by the simple well-known formula
\begin{equation}
    \sum_{n=1}^N n = 1+2+\cdots + N = \frac{N(N+1)}{2} \ .
    \label{sum r=1}
\end{equation}
On the other hand, while occasionally useful, the general formulas for the sums of higher powers are much less well-known. For example, one has
\begin{equation}
    \sum_{n=1}^N n^2 = \frac{1}{3} N^3 + \frac{1}{2} N^2 + \frac{1}{6} N \ .
    \label{sum r=2}
\end{equation}
In fact once a formula of the type of eq. \eqref{sum r=1} or \eqref{sum r=2} is given, its proof is a trivial exercise by using induction. However, \textit{guessing} the formula in the first place is non-trivial. 

Here we will provide a simple-to-remember method to recursively compute formulas for sums of a generic power $r \in \mathbb{N}$:
\begin{equation}
    \sum_{n=1}^N n^r \equiv S(N;r) \ .
    \label{sum r}
\end{equation}
In other words, our method allows us to obtain $S(N;r+1)$ from $S(N;r)$, and therefore recursively all the $S(N;r)$. In general, $S(N;r)$ will be a polynomial in $N$ of degree $r+1$.

\section{Faulhaber's Formula}\label{sec:faulhaber}

In the early 18th century, Jacob Bernoulli obtained a general formula for the sums of powers, now known as Faulhaber's formula,
\begin{equation}
    \sum_{n=1}^N n^r = \frac{1}{r+1}\sum_{j=0}^r (-1)^j {r+1 \choose j} B_j N^{r+1-j} \ ,
    \label{faulhaber}
\end{equation}
where the $B_j$ are the Bernoulli numbers with $B_1 = -1/2$. The proof of eq. \eqref{faulhaber} is straightforward starting from one of the forms of the Euler-Maclaurin identity \cite{spivey}.
While technically this is a solution to the problem, it is very hard to remember. Moreover, despite the apparent simplicity of eq. \eqref{faulhaber}, the computation of the Bernoulli numbers is itself non-trivial. Formulas are either implicit or involve double summations of complicated summands. For these reasons we look for an easier-to-remember method that one can perform with pen and paper without having to look up additional references. 

\section{A simple mnemonic for sums of powers}

Here we illustrate the mnemonic to recursively compute sums of powers, while in the next section we give a proof of the correctness of the method. Let's consider the concrete example of the computation of $\sum_{n=1}^N n^2$. Here we first do something \say{illegal}, that is we differentiate the whole sum with respect to $N$, and by this we mean that we differentiate term-by-term. Doing this we obtain the sum of one power lower $\sum_{n=1}^N n^2$, which we know:
\begin{equation}
    \frac{d}{dN} \sum_{n=1}^N n^2 \overset{?}{=} \sum_{n=1}^N \frac{d}{dn} (n^2) =2\sum_{n=1}^N n = N(N+1) = N^2+N \ .
    \label{miao}
\end{equation}
This is clearly not a valid operation, however one can formally think of it as a map between polynomials that happens to be given by a term-by-term derivative. In order to recover the original sum, one may exploit the analogy with differentiation and integration, whereby integrating reverses the derivative. Thus integrating the end result of \eqref{miao}, one obtains
\begin{equation}
    \int_0^N N^2 +N \,\, dN =\frac{1}{3} N^3 + \frac{1}{2} N^2 \ .
\end{equation}
Comparing with eq. \eqref{sum r=2}, this is almost the correct result, lacking only the linear term $N/6$. The procedure thus works as follows:
\begin{itemize}
  \item Differentiate the sum term-by-term and substitute the formula for the lower sum of power;
  \item Integrate the resulting polynomial;
  \item Add a linear term $CN$ and fix the constant $C$ by requiring that the sum give $1$ for $N=1$ as it should. 
\end{itemize}
Thus one is able to recursively obtain all the sums of powers for any $r$ starting from either $r=1$ or trivially $r=0$. In the next section we prove that this method indeed gives the correct result. This method can be easily remembered (\say{differentiate, then integrate and add a linear term}).

Now we compute the first few sums of powers as an example. Starting with $r=0$,
\begin{equation}
    \sum_{n=1}^N n^r=\sum_{n=1}^N 1 = N 
\end{equation}
trivially. Now we apply our procedure to the $r=1$ case. First we integrate and sum the resulting sequence:
\begin{equation}
    \frac{d}{dN} \sum_{n=1}^N n = \sum_{n=1}^N 1 = N \ .
\end{equation}
Then we integrate and add a linear term, finding
\begin{equation}
    \sum_{n=1}^N n = \frac{1}{2} N^2 + CN \ .
\end{equation}
Substituting $N=1$ we find $1/2 + C = 1$, that is $C=1/2$, so
\begin{equation}
    \sum_{n=1}^N n = \frac{1}{2} N(N + 1) \ ,
\end{equation}
which is the correct formula. Now for $r=2$ we have
\begin{equation}
    \frac{d}{dN} \sum_{n=1}^N n^2 = 2\sum_{n=1}^N n = N^2+N \ .
\end{equation}
Integrating and adding a linear term
\begin{equation}
    \sum_{n=1}^N n^2 = \frac{1}{3} N^3+\frac{1}{2} N^2 + CN \ .
\end{equation}
Substituting $N=1$ we have $1/3+1/2 + C = 1$, that is $C=1/6$, so that
\begin{equation}
    \sum_{n=1}^N n^2 = \frac{1}{3} N^3+\frac{1}{2} N^2 + \frac{1}{6}N \ .
\end{equation}
For $r=3$, we have
\begin{equation}
    \frac{d}{dN} \sum_{n=1}^N n^3 = 3\sum_{n=1}^N n^2 = N^3+\frac{3}{2} N^2 + \frac{1}{2}N \ .
\end{equation}
Again integrating and adding a linear term
\begin{equation}
    \sum_{n=1}^N n^3 =\frac{1}{4}N^4+\frac{1}{2} N^3 + \frac{1}{4}N^2 + CN \ .
\end{equation}
For $N=1$ we have $1/4+1/2 + 1/4 + C = 1$, that is $C=0$, so
\begin{equation}
    \sum_{n=1}^N n^3 =\frac{1}{4}N^4+\frac{1}{2} N^3 + \frac{1}{4}N^2 \ .
\end{equation}
For $r=4$,
\begin{equation}
    \frac{d}{dN} \sum_{n=1}^N n^4 = 4\sum_{n=1}^N n^3 = N^4+2 N^3 + N^2 \ .
\end{equation}
Therefore
\begin{equation}
    \sum_{n=1}^N n^4 =\frac{1}{5}N^5+\frac{1}{2} N^4 + \frac{1}{3}N^3 + CN \ ,
\end{equation}
where $C=-1/30$, so that
\begin{equation}
    \sum_{n=1}^N n^4 =\frac{1}{5}N^5+\frac{1}{2} N^4 + \frac{1}{3}N^3 - \frac{1}{30}N \ .
\end{equation}
In principle one can keep going and compute the formula for an arbitrary $r$.

\section{Proof of correctness}

Here we prove that the method explained in the previous section is correct. To do so, define
\begin{equation}
    S(N;r) = \sum_{n=1}^N n^r \ ,
\end{equation}
where $S(N;r)$ is a polynomial in $N$. Then the statement of the previous method is essentially that
\begin{equation}
    S(N;r) = CN + r \int_0^N S(N;r-1)\,dN 
\end{equation}
for some constant $C$. To show that this is true, we compute the integral using Faulhaber's formula:
\begin{align*}
    r\int_0^N S(N;r-1)\,dN &=r\int_0^N \frac{1}{r}\sum_{j=0}^{r-1} (-1)^j {r \choose j} B_j N^{r-j}\,dN =\\
    &=\sum_{j=0}^{r-1} (-1)^j {r \choose j} B_j \frac{1}{r+1-j} N^{r+1-j}=\\
    &=\frac{1}{r+1}\sum_{j=0}^{r-1} (-1)^j {r+1 \choose j} B_j N^{r+1-j}=\\
    &=S(N;r) - \frac{1}{r+1} (-1)^r {r+1 \choose r} B_r N =\\
    &=S(N;r) -(-1)^r B_r N
\end{align*}
That is
\begin{equation}
    S(N;r) = (-1)^r B_r N + r \int_0^N S(N;r-1)\,dN \ .
\end{equation}
This is precisely the statement that we meant to prove, with an explicit value for the constant $C$. While this may be useful, we prefer to determine $C$ by setting $S(N;r)=1$ for $N=1$ as it is easier to remember.

\end{document}